**ORIGINAL PAPER**

# Neyman-Pearson Hypothesis Testing, Epistemic Reliability and Pragmatic Value-Laden Asymmetric Error Risks

Adam P. Kubiak[1] · Paweł Kawalec[1] · Adam Kiersztyn[2]

**Abstract**
We show that if among the tested hypotheses the number of true hypotheses is not equal to the number of false hypotheses, then Neyman-Pearson theory of testing hypotheses does not warrant minimal epistemic reliability (the feature of driving to true conclusions more often than to false ones). We also argue that N-P does not protect from the possible negative effects of the pragmatic value-laden unequal setting of error probabilities on N-P's epistemic reliability. Most importantly, we argue that in the case of a negative impact no methodological adjustment is available to neutralize it, so in such cases the discussed pragmatic-value-ladenness of N-P inevitably compromises the goal of attaining truth.

## 1 Introduction

Many people are inclined to believe that the scientific method should be a reliable tool for reaching the aim of maximizing truth and minimizing falsity in a body of assertions. J. Neyman and E. Pearson's (see e.g. Neyman 1952) conception of hypothesis testing (N-P hereafter) may appear to be designed to address this goal by demonstrating the reliability of research methods to reach the truth in the long run. This supposed reliability is captured by error probabilities that reflect how often, given an iterated use of the method, the research process would end up with a false assertion. Nevertheless, the epistemic adequacy of N-P account is challenged by several philosophical-methodological problems raised by both philosophers (e.g.

✉ Adam P. Kubiak
  adampkubiak@gmail.com

1  Faculty of Philosophy, John Paul II Catholic University of Lublin, Aleje Racławickie 14, 20-950 Lublin, Poland

2  Department of Computer Science, Lublin University of Technology, Nadbystrzycka 36B, 20-618 Lublin, Poland

Sprenger 2013) and scientists (e.g. Ioannidis 2005). These problems undermine the possible attribution of truth-conduciveness of N-P and by that its scientificness. Not without reason has the method been conceptually interpreted by its designers as a decision-theoretic tool (e.g. Neyman 1957), which means a tool for drawing practical conclusions (decisions) from a mixture of epistemic and non-epistemic (social, cultural, economic, ethical, etc.) premises. Nonetheless, as long as N-P remains to be used as an element of the scientific method, the expectation of epistemic merit may remain valid. Perhaps the existence of the explicitly present pragmatic (social, economic, ethical, etc.) preferences for avoiding one type of error more than another in the N-P testing procedure (pragmatic value-ladenness) (see Sect. 2) makes N-P a more socially responsive method. But the question remains if it is worth the epistemic loss which stems, for example, from some kind of pragmatically motivated bias. In this paper, we examine whether N-P can be seen as principally satisfying, in a minimal sense, some general epistemic standards and how pragmatic value-laden uneven setting of error probabilities can influence it.

One may repel the problems of N-P's questionable epistemic value by simply treating N-P as a decision-theoretic framework, like Neyman himself vehemently did. Another means of defense would be to abstract from the pragmatic aspect of N-P and defend it as being truth-directed by modifying it and granting it a special philosophical interpretation (see e.g. Mayo, Spanos 2010). A different approach to defending the epistemic viability of N-P might be to admit that the epistemic reliability (ER hereafter) of the method is not the same as the epistemic interpretation of its outcome (Rochefort-Maranda, 2013). Additionally, the ER of N-P can be analyzed in the context of the relation between the pragmatic and the epistemic aspects of N-P. A combination of the last two ways is the route we follow in this paper.

A general reassessment of N-P's epistemic credibility is partially motivated by the inadequacy of the existing epistemic evaluations of frequentist hypothesis testing when applied to N-P. This inadequacy may be caused by the fact that the most common procedure used for testing data nowadays is an amalgamation of N-P and Fisher's procedures (see Peregonzales 2015). Consider the most general and basic analysis of ER of frequentist hypothesis testing conducted by greatly cited Ioannidis' (2005) work. His evaluation counts mistakes only among rejected nulls and therefore is more asymmetric than the original N-P, which also justifies accepting nulls[1]. Moreover, Ioannidis evaluates the performance of hypothesis testing in the context of manifold biases that accompany a practical use of statistical methods. These biases are not inherent elements of these methods, but are, for example, effects of the wrong use of the methods. Thirdly, he does not scrutinize the uneven importance of different types and sizes of errors from the epistemic perspective. Lastly,

---

[1] In the asymmetric approach, a rejection (with known error risk) can go only in one direction—that of rejecting $H$, but no such thing as a rejection of an alternative hypothesis and acceptance of $H$ is validated. N-P also includes rejection of an alternative and acceptance of $H$ with predefined error probability. This aspect should not be conflated with the aspect of asymmetric avoidance of errors of different type.

Ioannidis does not take into account the specific influence of pragmatic preferences on N-P's reliability, which is an inherent element of the method in the case of N-P.

Issues mentioned above indicate that there is an urgent need to re-examine the ER of N-P, which would: (I) be based on an adequate criterion for measuring ER, (II) measure the nominal reliability of N-P without taking into account biases that are external to the method, (III) evaluate the original version of N-P, not a modified or misinterpreted one, and (IV) account for N-P's pragmatic value-ladenness in the form of the unequal setting of error probabilities. Considering the fact that there is no such analysis in the literature, its outcome would provide a new—more basic and more complete—basis for discussions on the ER of N-P. We aim to conduct such an analysis of the credibility of N-P in tracking the truth that satisfies these four demands.

First, we review some rudimentary facts about N-P (Section 2). Then, we clarify the notion of ER and the intended scope of our assessment (Section 3). Next, we put forth the methodological assumptions of our analysis (Section 4) and provide an analysis of the ER of N-P (Section 5). Subsequently, we discuss the problem of uneven pre-study odds of true alternative hypotheses to true hypotheses tested (Section 6) and analyze the dependence of ER of N-P on the pragmatic values implemented in the form of the asymmetric setting of error probabilities (Section 7). Section 8 concludes by offering comments from a broader methodological-philosophical perspective.

## 2 Neyman-Pearson Hypothesis Testing

While we assume that the reader is already familiar with the contents of N-P, we overview some rudimentary facts with an emphasis on the uneven pragmatic importance of two types of errors assumed by N-P.

The application of a statistical test may result in four possible situations, two of which are unsatisfactory: (a) $H$ is true whereas the action taken is $B$, or (b) its complement $H^c$ is true, while the action is $A$ (see Neyman 1950, 261) (Table 1).

Table 1 shows that there are two kinds of random errors[2] associated with the two types of unsatisfactory situations:

    a. error of the Ist type $\alpha = P(rejectH|histrue)$, and
    b. error of the IInd type $\beta = P\big(acceptH|h'istrue\big)$[3]

where $h$ is a simple hypothesis being a particular instance of $H$ or being equivalent to $H$ (simple hypothesis $h'$ being a particular instance of $H^C$ or identical with it), and $1 - \beta$ is the probability that, given $h'$ [4], the sample point will fall in the rejection region specified for $h$[5]. This means the probability that the test will detect the falsehood of $H$ (reject it) when the true hypothesis is c—is what Neyman called the *power* of a test (Neyman 1952a, 55; Neyman 1950, 267–268)[6]. Power, as a function of the point hypothesis from possible hypothesis space, is the essential category used for assessing which test to choose.

---

[2] The verdict on taking a particular action is random because it depends on the random variable(s) determining the position of the sample point. Due to this, there is no inconsistency in considering the probability of the verdict having a certain property, such as being erroneous (Neyman 1950, 56-57).

[3] Today it is standard to use "$\beta$" to represent the probability of making an error of the $II^{nd}$ type but in N-P's original notation "$1 - \beta$" was used to denote it while "$\beta$" denoted power. We use contemporary notation to enhance readability and ease-of-comparison with other work in the area.

It is important to distinguish statistical hypotheses, which are mathematical models, from physical hypotheses, which are propositions belonging to some field of science. A physical hypothesis is a statement about physical reality as seen from the perspective of a definite scientific discipline (and its language) within which testing is performed[7]. The statistical hypothesis is a statement expressed in statistical terms that are assumed to represent mathematically the physical hypothesis in question. An example is the statement that a die is fair, as compared to the statement that the model of sampling probability distribution under experimental design used to verify the statement has such and such parameter(s) value(s). The latter is meant to represent the former, but the same statistical hypothesis—definite model of probability (density) distribution—may stand for different physical hypotheses belonging to different fields of science[8]. Still, a definite probability (density) distribution has the same mathematical properties regardless of what it denotes in the physical world. This means that it can be analyzed regardless of its semantic interpretation related to the physical context of a particular test. The same holds for the particular types of error: statistical laws can be applied to analyze them regardless of the physical contexts of a particular testing situations to which the errors apply. These two facts are important for our analyses.

Keeping the value of $\alpha$ error at the desired nominal (theoretical) level is unproblematic, as the researcher decides before the research under what significance level $\alpha$ a test procedure will be executed. The $\beta$ nominal error probability depends in turn on a fixed instance of the alternative hypothesis, the distribution of a test statistic—which is determined by values of sample size and population variance—and obviously on the chosen value of $\alpha$ that determines the rejection region. One important consequence of that is, at least in N-P, that, all other things being equal, an increase of $\alpha$ results in a decrease of $\beta$ for any particular parameter value belonging to $H^C$.

---

[4] Symbols $H$ and $H^C$ refer to sets of values (of sampling distribution model's parameter(s)), $h$, and $h'$ to particular values of it. Keep in mind that the meaning of $H$ and of $h$ is not limited the so-called "no effect" hypothesis, although, dependent on the case, they may denote this type of statement.

[5] The presented concept of errors is the most general approach that covers various test cases (e.g. with different distributions) and interval estimation. It does not cover errors of a different nature, like measurement instrument's errors, biases, or model assumption errors (like false normality or independence assumption).

[6] The value of power refers to departure from $H$ that is today understood also as the "expected minimum effect size" (Peregonzales 2015); this understanding is a consequence of the fashion of identifying $H$ with the substantial hypothesis of "no effect", "no difference", etc., which is not necessary as far as N-P is considered. The cited phrase means: the value of the alternative hypothesis' parameter that represents minimal value of the departure from the hypothesis tested that is of the researcher's interest (assuming the test statistic is ubiased).

[7] In the literature, it is also called "scientific hypothesis" (Hurlbert, Lombardi 2009, 313)

[8] Whether we speak of the statistical hypothesis, or the physical hypothesis (assumed to be represented by an adequate statistical hypothesis) should be clear from the context.

**Table 1.** Four possible outcomes of a statistical test

| True Hypothesis | $H$ | $H^C$ |
|---|---|---|
| Action taken | Description of the situation | |
| $A$: accept $H$ | Satisfactory | Error |
| $B$: reject $H$ | Error | Satisfactory |

Thus these two types of error probabilities will not necessarily have equal values and the difference will reflect a pragmatic discrimination of the importance of each type. Neyman considered this fact as consistent with the realm of application, where the importance of avoiding these two types of error is strikingly unequal:

"The adoption of hypothesis $H$ when it is false is an error qualitatively different from the error consisting of rejecting $H$ when it is true. This distinction is very important because, with rare exceptions, the importance of the two errors is different, and this difference must be taken into consideration when selecting the appropriate test" (Neyman 1950, 261).

## 3 The Epistemic Goal and the Intended Scope of our Assessment

Our understanding of the goal of attaining the truth is quite general: we consider it as the end-point of the epistemic process of scientific cognition. Although this goal may pose an unattainable ideal (see e.g. Rescher 1999), many philosophers would agree that scientific claims are rational because a method of scientific justification reliably gets us closer to the ideal situation. The *truth goal* as cast in terms of ER of a test procedures and measured by the frequency of true acceptances/ rejections, would thus mean that:

For all $H$ s considered by $S$ (individual scientist or a community) during a specific sequence of empirical studies in which N-P is used to perform a test, if $H$ is true, then $S$ accepts that $H$, and $S$ accepts that $H$ only if $H$ is true.[9]

N-P's characteristics may be insufficient to ensure that the actual error probabilities are satisfactorily close to the expected outcomes of a reliable method—the extent to which this reliability is close to the ideal of the above-defined truth goal. This means, for example, that a small $\alpha$-error-rate rate may not in fact yield a small error probability, if disturbing biases from outside the method (e.g. a bias of the measurement tool) are present. That is why we think of N-P's error probabilities as *nominal* error probabilities. In this paper, we concentrate on the investigation of the nominal ER. We assume that intuitive, absolutely minimal nominal epistemic requirement for a method of statistical inference would be that it will not lead to false assertions more often than to true ones, as well as that its ER will not remain at the level of a random guess. The difference between actual and nominal ER reveals that the critical examination of the ER of N-P can fall into one of the following general categories. First, (i) it can be an investigation of the epistemic adequacy of N-P procedures as such (e.g. Jaynes 2003). (ii) It may focus on phenomena and/or circumstances that affect the use of the method, for instance the effect of the publication bias (e.g. Dickersin 1990) or researchers' ignorance of the method and/or its (mis-)application (e.g. Gigerenzer, Merewski 2015). In general, this kind of analysis examines socio-psychological factors that are a non-inherent part of the N-P itself. (iii) Such an examination can concern the methodological insufficiency of the N-P. It may include, for example, issues of the need for complementary statistical tools that address the emerging new specific research circumstances, the need for some additional tools for amalgamating outcomes, or the issue of formulating the semantic content (physical denotation/interpretation) of considered statistical hypotheses, for which individual disciplines are responsible.

---

[9] Obviously, this statement is categorical but statistical inference always has an uncertain, probabilistic element therein. The truth goal is an ideal, to which closeness can be assessed by measuring a statistical method's efficacy. See also David (2001).

We wish to stress that in our examination we do not intend to touch upon the problems related to (ii) and (iii). We believe the topic we address here is more rudimentary and can be investigated independently from (ii) and (iii). Obviously, the nominal epistemic reliability of N-P is a necessary, but insufficient condition for a successful realization of the epistemic goal. Several other conditions belonging to (ii) and (iii) must be met as well, but they should not be treated as stemming directly from the inherent features of the method as such. By concentrating on (i) we, for example, do not consider problems of the amalgamation of outcomes related, for instance, to the socio-psychological file drawer effect (cf. Rosenthal 1979), a type of bias where "negative outcomes" (nulls) remain unpublished. We also do not aim to solve the problem of formulating the semantic content of hypotheses (within this or that discipline), which would be consistent with the physical theory in force (issue of type iii). Certainly, the influence of pragmatic values, in general, also embraces the issue of the influence on the semantic content of hypotheses/theories. One of the most striking examples is how Lysenkoism (natural sciences subjected to Marxist-Leninist philosophy) impeded the development of biological disciplines in the Soviet Union (see e.g. Soyfer 1994). This and similar circumstances seem to suggest that social values may interfere with scientific research by promoting false hypotheses and false theories. Yet, the indicated cases of societal influences apply to the semantic content of admissible, or accepted, hypotheses—they refer to the empirical interpretation of mathematical statements. By contrast, the type of influence considered here is the influence of values on the epistemic performance of a procedure as captured by error probabilities linked to different types of error.

No less important clarification of our goal is to specify what kind of examination of (i)-type (evaluation of epistemic adequacy of N-P as such) we do not follow in this paper, and to narrow the scope of the possible (i)-type topics that we concentrate on. We do not intend to evaluate N-P on particular evidence, as N-P is essentially concentrated on pre-trial reliability (see Graves 1978, 6-7). This does not mean that this evaluative perspective could not be adopted. For example, Spielman (1973) defines reliability indexes as the posterior probability of a correct action given a definite outcome. A more recent example of such an epistemic interpretation of the single outcome of the N-P procedure is exposed by Deborah Mayo (1996) accompanied by Aris Spanos (e.g. Mayo, Spanos 2006), who presented the most current epistemic reading of N-P and philosophical defense of the frequentist statistical paradigm (see Mayo 2018). Following the argumentation of Rochefort-Maranda (2013), we never­theless stress that an epistemic interpretation of N-P's single outcome (a measure of evidential support for a particular case of application) differs from an assessment of the epistemic credibility of the procedure. The strength of evidence for or against a particular hypothesis is a result of the quality of the test and its output. According to Rochefort-Maranda, "[o]ne way to distinguish both concepts (level of evidential support and level of its credibility) is to realize that the credibility of the support does not depend on the actual output of the instrument whereas the degree of sup­port does" (Rochefort-Maranda 2013, 11). Rather than solve what could be a post­experimental epistemic interpretation of a particular outcome we seek to analyze N-P's general epistemic credibility.

As for narrowing down the scope of our analysis, we do not intend to examine any of the existing (i)-type issues relating to an analysis of N-P in comparison to other frequentist approaches, like Fisher's (1956) asymmetric frequentism, or Bayes­ian approaches (see Romeijn 2017). We intend to take a closer look at N-P itself, without deciding which of the approaches is better in any respect.

## 4 Methodological Assumptions of the Analysis

Given Neyman's view on N–P as a theory of how one should make decisions and not a theory of how one should change beliefs (see e.g. Neyman 1957) the ER of N-P may seem irrelevant. Pragmatic evaluation is based on pragmatic purposes and the methodology of N-P testing itself: a well-defined N-P test has determinate error probabilities, and the question of whether a given N-P test is sufficiently useful for a given pragmatic purpose is easily answered by consulting its error probabilities and comparing them to the criteria of sufficiency that the purpose at hand warrants. Nev­ertheless, if an epistemic evaluation of N-P's performance is viable then it require the perspective from outside, not inside, the method. A classic example of an outside evaluation of the ER of frequentist hypothesis testing[10] is Ioannidis' (2005) assess­ment. It took into account an undefined within the frequentist framework pre-study probability of the hypothesis being true and used the concepts of positive predictive value and pre-study bias to assess the performance of the testing procedure.

At the outset, we use the same criterion for the analysis of the epistemic perfor­mance of N-P as Ioannidis, namely the probability that the alternative hypothesis would turn out to be true once it has been accepted with statistical significance; this is called its "positive predictive value"[11] (*PPV*) (2005, 696). PPV is dependent on the values of both types of error, as well as on the ratio $R$ of alternative hypotheses that are true to the alternative ones which are false among those investigated in the respective field of science[12]:

$$PPV = P\left(H^C | accept H^C\right) = (1 - \beta)R/((1 - \beta)R + \alpha). \tag{1}$$

---

[10] Although not exactly of N-P, as noted in the Introduction.

[11] Instead of using the term "alternative", Ioannidis spoke of the hypothesis of the existence of an effect in contrast to the hypothesis of no effect, which is colloquially called the "null" hypothesis.

$PPV$ is a concept used in medical statistics and its counterpart is the concept of "negative predictive value" ($NPV$), i.e. the probability that the tested hypothesis would turn out to be true once it has been accepted with statistical significance (see Altman, Bland 1994):

$$NPV = P(H | accept H) = (1 - \alpha)/((1 - \alpha) + \beta R). \tag{2}$$

while Ioannidis uses only $PPV$, we are going to use $PPV$ and $NPV$ because of the symmetry N-P is in accepting hypotheses.[13]

We assess N-P's credibility from the perspective of $PPV$ and $NPV$, assuming Neyman and Pearson's important requirement that

(A1) a test should be at least *unbiased*.

A test is unbiased when its power against any alternative point hypothesis is at least as high as the type I error (see Neyman, Pearson 1936, 210–211). The basic rationale of an unbiased test is that it avoids cases where the probability of accepting $h$ when it is false (and the true hypothesis is $h'$) would be higher than the probability of accepting it when it is true (see Neyman, Pearson 1936, 210–211).

The case of $1 - \beta = \alpha$ can be plausibly identified as a value corresponding to the case of probability of detecting infinitely small departure from $H$ (or committing infinitely small error of the IInd type). Power $(1 - \beta)$ is the probability of obtaining an outcome that falls under the $\alpha$ rejection interval given the assumption that the unknown quantity has one definite value belonging to the admissible hypothesis' parameter space. In the case of an unbiased test, the closer this true value is to $H$, the lower is the probability of observing an outcome that falls under the $\alpha$ interval. This probability is lowest when it is equal to $\alpha$, which is the case if the true value is identical to the value of $H$ for which $\alpha$ rejection interval was set, or if the departure of the value of $H$ (for which $\alpha$ rejection interval was set) from the truth is infinitely small. In sum, considering possible situations of a value of power to detect the discrepancy of an accepted $H$ from the true value that is an instance of $H^C$, the case of the lowest possible power (equal to $\alpha$) can be plausibly understood as standing for an infinitely small error of false acceptance of $H$. Below we argue that it is not meaningful epistemically to consider this case in the evaluation of the ER of N–P.

If standards for avoiding errors are very stringent the goal of attaining the truth is promoted, but at the same time it is blocked at a more fundamental level (Steel 2010): very exacting standards suspend conclusions; devoting more and more resources to make the judgment more and more accurate not only suspends making

---

[12] The pre-study probability that the alternative is true ($P\left(H^C\right)$) is therefore $R/(R + 1)$ (see Ioannidis 2005, 696).

[13] See footnote 2.

that judgment but also deprives the researcher of potential cognitive resources to be used for other research questions. Because our cognitive resources are limited, a trade-off arises between the need to avoid mistakes and the need to be able to effectively scrutinize hypotheses in finite time by limited resources. Setting the threshold for being wrong but acceptably close to the truth allows the researcher to continue testing other hypotheses; differentiated error probabilities reflect such a trade-off between the stringency of a test and the need to test the hypothesis in a finite time and concerning available resources in a given research context. For different research contexts, this balance can be different due to the differences in limiting resources connected to pragmatic contexts of particular research (see Steel 2010, 27–28). The conclusion of that argument is that minimizing errors exaggeratedly in a given test is epistemically unfavorable; in particular, making an effort to minimize errors to an infinitely small level to secure infinitely high accuracy of conclusions (which means to secure perfectly errorless/accurate conclusions) is epistemically irrational. Therefore, when the ER of N-P is evaluated[14] taking into account the possible cases of infinitely small errors of the second type is meaningless. As we argued in the previous paragraph, consideration of this type of error would mean consideration of $1 - \beta$ being equal to $\alpha$, thus

(C1) taking into account the case of $1 - \beta = \alpha$ in the valuation of N-P's ER is meaningless.

Cases of infinitely small errors can be treated as epistemically meaningless. The minimal reliability does not need to be achieved in these cases—*PPV* and *NPV* values for them can be ruled out as being epistemically irrelevant. Therefore, epistemically relevant cases to be investigated are those when $1 - \beta > \alpha$.

We wish to address two emerging issues concerning the use of *PPV* and *NPV* as indicators of N-P's ER. The first problem is the use—after Ioannidis[15]—of the concept of prior probability in our assessment of *PPV* and *NPV*. The use of prior probabilities is not an element of N-P. But to measure the method's ER one does not have to stick to the tools and concepts that are part of the method itself. By analogy, one cannot use a rocket's engine to measure this engine's noise level. Also, it is not a goal of an engine to make noise. But it does not mean that engine does not make noise measurable by a suitable tool. The use of priors in the assessment of ER of N-P does not assume the use of priors in the N-P inferential procedures itself, so frequentist arguments for N-P as such should have no force as a potential critique of our approach, which is taken from an outside perspective.

Perhaps for orthodox frequentists, like Neyman himself, our analysis may appear irrelevant or uninteresting due to their decision-theoretic interpretation of the goal of N–P and thus the sole interest in its decision-theoretic, not epistemic, reliability. But since an assessment of N-P from an epistemic, outside perspective is still possible, its outcome may be interesting for unorthodox frequentist philosophers of science and statisticians, as well as for Bayesians. The perspective we adopt is Bayesian-like because of reference to the pre-study probability of a hypothesis yet it is not Bayesian for two reasons. Firstly, from the Bayesian perspective, the hypothesis' probability could be different from case to case, provided one employs more specific information that allows assessing the prior probability of the hypothesis in a particular case. In our analysis, as in Ioannidis', a general method's performance is analyzed, making precise assumptions concerning the prior probability that refer to the context of individual research would only allow us to infer ceteris paribus conclusions for particular research cases. A minimal assumption about $R$—whether it is equal to, greater or smaller than 1—allows for sufficiently general statements about N-P's reliability, and at the same time to distinguish types of cases that are relevant for the goal of analyzing the influence of pragmatic value-ladenness. Therefore, in our analysis, we assume that the only—but also sufficient and adequate for our purposes—knowledge at hand is that $R$ for the type (e.g. discipline) to which the research belongs is equal, greater, or smaller than 1.

---

[14] Additionally, in the simplest case where the sampling distributions of $H$ and an alternative hypothesis are of the same type and have common variance, the test is uninformative as there is no separation between sampling distributions (they coincide) when $1 - \beta = \alpha$.

[15] What we already signalled, Ioannidis derived pre-study probability that the alternative is true from $R$.

In the light of this type of scarce information, all that can be said (and is adequate from the perspective of the purpose of our analysis) is that the pre-study probability of the hypothesis tested during the research of the given type is equal to, smaller or greater than 0.5. Ioannidis considered cases of $R$ smaller, or bigger than 1 dependent on scientific discipline considered. This is covered by our analysis. Yet, one may also want to consider N-P's reliability in the broadest possible sense—as applied in manifold research situations. Following Ioannidis, let us assume that for some branches, or types of studies, $R$ is smaller, and for some other, greater, than 1. The reference class of disciplines in such a case is so broad that one has no reason to think that most possible hypotheses tested are true, nor that most of them are false. We assume, that in such a case the most plausible move would be to refer to the principle of indifference and assume that the pre-study odds are equal to 1, which is the third case covered by our analysis.

The second difference from Bayesian approach is that in our analysis we do not refer to the concept of probability of the hypothesis given definite evidence: *PPV* and *NPV* do not inform about the probability of a hypothesis being true given the evidence obtained but given the fact of acceptance of the hypothesis. Following the interpretative assumptions as described above we will use $R$ and $P(H)$ as admissibly equivalent measures of testing conditions.

The second problem concerning the use of *PPV* and *NPV* as indicators of N-P's ER that we wish to address is the following. A definite predictive value, for example *PPV*, always refers to a point hypothesis that corresponds to definite values of error probabilities of $\alpha$ and $\beta$ so the alternative hypothesis, to which *PPV* refers, is not a composite, but a point alternative hypothesis that corresponds to the value of $\beta$. Hence is of a definite distance from the point hypothesis that corresponds to the value of $\alpha$. That is why, strictly speaking, the probabilities used to calculate *PPV* and *NPV* only refer to those two particular possible values of point hypotheses that are the instances of $H$ and $H^C$.

In our analysis, we assume that the prior probabilities (on which we operate when talking about definite predictive values) refer not only to parameter values for which values of $\alpha$ and $\beta$ stand but to them taken jointly with more extreme parameter values, that is to those values that are farther from the hypothetical falsely accepted value in the case of a given type of error. Owing to this, the situation is as follows: when speaking of the value of $NPV$, we are thinking of a statement about the predictive value under the assumption that a true value of the parameter is within the bounds of a range that corresponds to error probability equal to $\alpha$ or error probabilities smaller than $\alpha$. Similarly, in case of the value of $NPV$ we speak of the true value of the parameter being the one that corresponds to an error probability that equals $\beta$ or to smaller error probabilities. Additionally, we assume that for all the more extreme parameter values error probabilities $\alpha$ (and $\beta$, respectively) that correspond to these values are equal to those particular values of $\alpha$ and $\beta$ that are considered to be the basis for calculating definite values of $PPV$ and $NPV$.[16] Error probabilities for these more extreme parameter values are, by standard, smaller, but for computational simplicity such a simplifying assumption may suffice. Under this assumption, particular values of $\alpha$ and $\beta$ are the upper bounds of nominal error probabilities for the considered possible ranges of true values of the unknown parameter, assumed under the stated values of $\alpha$ and $\beta$. So the predictive values ($PPV$ and $NPV$) calculated from those stated error probabilities and some stated value of $R$ are the lower bounds for the predictive values.

## 5 The Epistemic Reliability of N–P

Ioannidis was particularly worried about the low value of $PPV$ in the case of some disciplines that are, according to him, characterized by a low $R$ and underpowered study designs. Consider $\alpha = 0.01$, very low power $1 - \beta = 0.02$ (N-P's unbiasedness condition is satisfied), and very small pre-study odds $R = 0.02$. In such a case $PPV \cong 0.04$, so it is indeed very low, but at the same time $NPV \cong 0.98$, which is very high. As a consequence, if such a testing condition is present in an iterated use of N-P, one will not commit errors very often. The fraction of true acceptances among all the acceptances will be approximately equal to 0.97, hence quite satisfactory despite the very low $PPV$. If power was higher, the overall error rate could only be smaller.

Cases like the ones above cannot be properly captured if one separately uses $NPV$ and $PPV$ as indicators of epistemic credibility. To be able to grasp the described case and thus evaluate the reliability of N-P itself, rather than a hybrid of N-P with the Fisherian approach, one has to consider $NPV$ and $PPV$ jointly[17]. It is possible by looking at the value of the total probability of accepting a true statement ($P_t$). It equals the sum of the probabilities of two mutually exclusive events. first, the case of $H$ being true and $H$ being accepted, and, second, the case of $H^C$ being true and $H^C$

---

[16] This is analogical to Ioannidis' simplifying assumption of equal power to find all the true effects existing in the field of study.

[17] It is consistent with the symmetric approach of N-P to treat both rejection and acceptance of $H$ as a research outcome.

being accepted. This concept, by the definition of conditional probability, (1) and (2), yields[18]

$$P_t = P(acceptH)NPV + P\left(acceptH^C\right)PPV. \tag{3}$$

By the definition of the conditional probability $P_t$ can be restated as

$$P_t = (1 - \alpha)P(H) + (1 - \beta)P\left(H^C\right). \tag{4}$$

Given N-P's symmetry $P_t$ appears to be the most intuitive, simple, and adequate measure of ER of N-P that is anchored in standard *PPV* and *NPV* measures. Its usage can be illustrated by a simple example. Imagine the situation of tests of hypotheses with $\alpha = 0.05$ wherein half of the cases the hypothesis tested is false ($P(H) = 0.5$) and true values are at some distance from the hypothesis tested not smaller than a certain value that marks $1 - \beta$, say, 0.1. Then $P_t$ will be equal to 0.525, which means that, for example, if 1000 hypotheses were tested, the nominally expected number of true research findings would be at least 525.

As postulated in Section 3, the minimal epistemic requirement for a method of statistical inference is that it will not lead to false assertions more often than to true ones, as well as that its ER will not remain at the level of a random guess. Therefore, the condition of minimal ER is that

$P_t > 0.5$ for any testing situation under the assumption of (A1) and condition (C1)

(5)

Given (A1) and (C1) the condition of minimal epistemic reliability (5) holds for any test when $R = 1$ ($P(H) = 0.5$). That is because if $(1 - \beta) > \alpha$ then $\{(1 - \alpha) + (1 - \beta)\}/2$ is greater than $1/2$. Alas, condition (5) is not satisfied when $R < 1$ ($P(H) > 0.5$). For example, if $P(H) = 0.95$, $\alpha = 0.55$ and $1 - \beta = 0.99$, then $P_t = 0.48$. In this example, the probability of error of the Ist type is greater than 0.5, but the maximal admissible level of $\alpha$ was not stipulated by Neyman and Pearson. Moreover, Neyman claimed that it may in some cases be treated as the less important type of error (Neyman 1971, 4)[19]. Condition (5) is not satisfied also when $R > 1$ ($P(H) < 0.5$). For example, if $P(H) = 0.33$, $\alpha = 0.05$ and $\beta = 0.94$, [20] then $P_t = 0.35$. The conclusion is that condition (5) holds only for $R = 1$ ($P(H) = 0.5$).

---

[18] The reader may keep in mind that *PPV* and *NPV* do not inform about the probability of a hypothesis given the evidence obtained, but given the acceptance of the hypothesis.

[19] What one might arguably assert as the principle that stipulates a lower bound for $\alpha$ in such a case is unbiasedness applied in a reversed form, namely $1 - \alpha$ must be at least as high as $\beta$; this condition is fulfilled by the example given.

[20] Note that this example represents the case of an underpowered test, but is, still, consistent with (A1) and (C1).

## 6 Asymmetry of Errors, Hypothesis' Probability, and Epistemic Reliability

Securing the desired level of the Ist type of error can be arbitrary, easier than securing the level of the IInd type of error[21]. Therefore, in practice $1 - \alpha$ is greater than $1 - \beta$, and then the greater $P(H)$ is, the greater is also $P_t$. It then appears that it would be better epistemically to test a hypothesis that is deemed to be true rather than false[22]. It turns out that the supplementation of N-P with a pre-experimental assumption about the probability of the tested hypothesis is supported by Neyman and Pearson's suggestion expressed in one of N-P's foundational papers. The suggestion they made is that a researcher usually has reasons for believing the tested hypothesis is true or that the true value differs from the point hypothesis stated "in certain directions only", which means that a definite composite hypothesis is true:

> "the hypothesis whose probability we wish to test is that $\Sigma$ is a random sample from $\Pi$" (Neyman, Pearson 1928, 178) and "[i]t is true that in practice when asking whether $\Sigma$ can have come from $\Pi$, we have usually certain *à priori* grounds for believing that this may be true, or if not so, for expecting that $\Pi'$ differs from $\Pi$ in certain directions only" (186).

The suggestion can be expressed as:

> (A2) A researcher usually has a priori grounds to believe that point or composite hypothesis being tested is true

This can be understood as a contextual assumption about the prior degree of belief in the hypothesis being tested. It can be translated into the probabilistic statement: $P(H)$ is higher than $P(H^C)$. Although the assumption is not implemented in the statistical inference scheme itself (N-P does not operate on the prior or posterior probability of a hypothesis), it can be understood as a precondition under which the method is used. Making this assumption does not lead to inconsistency with N-P. After all the assumption (A2) does not entail the technical application of its probabilistic statement in the method's procedures itself. Additionally, an evaluation that ignores the pre-assumed proper context of use would, for instance, resemble an assessment of the performance of binoculars at night, while it is known that they are presumed to be used in the presence of light. Hence, evaluation in the context of (A2) seems more adequate in the sense of taking into account Neyman and Pearson's concept as a whole, thus as a proposal of certain methods working under a certain assumed condition. In any case, even if (A2) was not stated as a normative circumstance (that must be met if one wants to use the method), it can be assumed to be conceived by Neyman and Pearson as the inherent part of the usual research practice.

---

[21] For example, $1 - \beta$ is dependent on the variance of the studied quantity in the population.

[22] In scientific practice the tested hypothesis $H$ is often associated with the "null" hypothesis, which states that the size of an effect being investigated is equal to 0. But, when scientists perform regression analyses, for example, they usually suspect that it is not equal to 0. In such cases, it is possible to perform data transformation so that the hypothesis they suspect to be true becomes the "null" hypothesis. For example, there are usually maximal and minimal values possible to be observed within a given experimental scheme. One can mathematically reverse these values. This will change the way particular outcomes are labeled by values of the random variable and therefore change the way the tested statistical hypothesis is defined, but will not change the physical meaning of the redefined statistical hypothesis.

Whether this assumption is representative of what happens in practice falls into the already distinguished set (ii) (phenomena and/or circumstances that affect the use of the method; see Section 3) of issues that we do not intend to answer. But let us briefly address this question. It seems fairly natural to admit that in general, when one comes up with a physical hypothesis to be investigated, one suspects that the state of affairs postulated might be true rather than false. This can be realized in two ways: either the hypothesis is a consequence of a theory that has been accepted so far as the most compelling, therefore rather true than false, or the researcher's background knowledge prompts her to suspect that some epistemically interesting fact might be the case, which pushes them toward formulating a new hypothesis. Some scientists might then make the effort to transform variables so that the intuitively true hypothesis becomes the alternative hypothesis: because the hypothesis tested then becomes—so to speak—counterfactual in the researcher's eyes, there is a greater hope that the observation will yield a substantial result. If, on the other hand, the researcher expected the tested hypothesis to be true, they would expect the observation to speak in favor of the tested hypothesis, which would mean she would expect a result difficult to be published (see Rosenthal 1979); this, of course, assumes an asymmetric approach, which is not the case in N-P.

Regardless of whether it can be affirmed that assumption (A2) is met in practice, (A2) can be regarded as a normative requirement that specifies the basic state of affairs that constitutes the presumption about the typical context of the application of N-P. This presumption can be satisfied using N-P. Firstly, a researcher can indicate, before the research, whether the hypothesis under investigation is expected to be true or rather false[23], without providing a more precise probability statement and thus without including it in the statistical inference. Secondly, even if the alternative is more probable, it is in principle possible to mathematically redefine the alternative as the hypothesis to be tested through a transformation of the test or data, and without any loss of empirical adequacy. In particular, it is easy to do that in the case of directional tests of a composite hypothesis. Nevertheless, in the case of test of a point hypothesis, when the set of admissible hypotheses is continuous, the only reasonable precondition appears to be $R \geq 1$: to assume that the prior probability of $H$ is greater than the prior probability of its complement seems mathematically incorrect in such a case.

---

[23] This is what happens in scientific practice. Usually, some formerly justified theory existing within a discipline, or an expert's intuition, foresee some hypotheses to be true, and the role of testing is to check whether the suspected state of affairs is indeed the case.

## 7 The Pragmatic Factor and Epistemic Reliability

From what we just have discussed it follows that pragmatically driven asymmetry of probabilities of two types of error might seem naturally epistemically beneficial under the condition (A2) assumed by Neyman and Pearson. It seems that from an epistemic perspective it would be better to select for testing a hypothesis deemed to be rather true. But at the same time, pragmatic asymmetry in avoiding the two types of error, along with the fact of the relative ease of controlling $\alpha$ (technically, its value can be freely set) as compared to $\beta$ justifies N-P's recommendation to design the tested hypothesis as one of the possible alternatives of which erroneous rejection would be pragmatically worse. That is why "[i]n an example of testing a medical risk, Neyman says he places 'a risk exists' as the test hypothesis since it is worse (for the consumer) to erroneously infer risk absence" (Mayo 2018, 341). This means that from the epistemic perspective it would be better if this pragmatically driven decision on which of the alternatives to place as the hypothesis tested coincided with the hypothesis that is regarded the more probable as compared to the alternative. Or, more generally, it would be better if the more probable physical hypothesis is at the same time the one of which false rejection should be more avoided than the false rejection of an alternative considered. The crucial question is then if this coincidence naturally occurs or can be stipulated. Below we argue that such coincidence is not sufficient for an increase of $P_t$ and that such coincidence is haphazard and cannot be stipulated.

In some cases, pragmatic value-laden asymmetric error risk setting can directly improve the ER of N-P. Assuming (4), (A1) and (A2), one might expect that, with resources fixed, an increase of $1 - \alpha$, at the cost of a decrease of $1 - \beta$, should improve the ER of N-P. $P_t$ increases if the modulus of the magnitude of the increase of $1 - \alpha$, multiplied by $P(H)$, is greater than the modulus of the magnitude of the decrease of $1 - \beta$, multiplied by $P(H^C)$, but such a case is not a rule. Holding sample size fixed the rise of $P_t$ would be guaranteed only if the loss of power was at most as big as the gain of $1 - \alpha$, but this will not typically be the case. For example, if the distribution of a sampling statistic is close to the normal distribution and the $1 - \alpha$ integral gets larger, it does so at a lower rate than the power integral diminishes. Because of that, for $P_t$ to rise P $(H)$ would have to be sufficiently high so that the modulus of the magnitude of the increase of $1 - \alpha$, multiplied by $P(H)$, is greater than the modulus of the magnitude of the decrease of $1 - \beta$ multiplied by $P(H^C)$. Therefore, the discussed coincidence is not sufficient for $P_t$ to rise. This means that the promotion of avoidance of erroneous rejections of a more probable hypothesis by a pragmatically driven asymmetry in error probability does not suffice to improve ER of N-P: $P(H)$ would have to be sufficiently high for the improvement to take place.

Even if the probability of a more probable hypothesis would be sufficiently high, there remains the earlier mentioned problem of lack of (the possibility of stipulating) the discussed convergence. Even if the more probable physical hypothesis is mathematically defined as the one that is tested ($H$) it will not necessarily coincide with it being the one of which false rejection should be more avoided than the false rejection of the alternative, from a pragmatic perspective. First of all, the concept of a more important error as such is not uniquely connected to errors of the first kind. Whether it will be the more probable hypothesis of which wrong rejection will be more important to avoid depends on the context of research. There may be cases for which one will prefer to avoid an error of the IInd type (Neyman 1971), while the hypothesis tested will be believed to be rather true. Neyman's simple example of a research context for which a IInd type of error would be more important was an investigation of quality control:

"Practical situation of this sort may be illustrated by the case of production in a factory. As long as the process of manufacture is characterized by the mean value $\xi \leq \xi_0$ of a certain characteristic $X$ of the product, the situation is satisfactory and no changes are necessary. On the other hand, if the mean value of $X$ becomes $\xi > \xi_0$, it is imperative to stop the process and to readjust the machines" (Neyman 1971, 4).

Suppose that the characteristic of interest is, for example, tolerance of a resistor (its departure from nominal values of resistance). Imagine $H$ to be the statement that the tolerance of resistors produced by a particular line does not, on average, exceed $\xi_0$. The manufacturer is expected to produce resistors of the given maximal tolerance as a part of a contract with a large company that needs them to produce laboratory devices in which precision is crucial. The contract specifies that $H$ must be satisfied. Otherwise, a very high fine is anticipated, which jointly with the loss of the market's trust in the manufacturer would almost surely ruin their whole business. If the manufacturer suspects $H$ to be false, they may withdraw the batch of resistors and readjust the machines, which is assumed to generate a reasonably low economic loss compared to selling a product that did not meet the standard. In such circumstances, the more important error to avoid is the error of falsely asserting $H$ when it is not the case. Having in mind that this error is the additive complement of the power to detect the falsity of $H$ (for a given $h'$), "(…) the desirable property of the test of $H$ is as high a power as practicable, perhaps with some neglect of the probability of rejecting $H$ when true" (Neyman 1971, 4).

At the same time, it seems that under the above type of research context the default assumption here is that $H$ is rather true. Otherwise, the manufacturer would not sustain the production of probably useless resistors in the first place. The idea of such research is that although the production process is in principle expected to yield usable resistors there is always a minor risk of something going wrong with it and because of this risk that iterative control is needed. Even if one made the alternative hypothesis to be the tested one and by that error of the first type to be more important, the redefined tested hypothesis could not be any more regarded as the more probable one.

A possible lack of the required correspondence can also occur in research cases different from quality control. For example consider the test of whether a newly created, improved version of a drug is non-toxic. In this type of research, the IInd type of error may be regarded as more important from the perspective of the role of outcomes of medical research for society. Simultaneously, if the previous version of the drug was confirmed non-toxic, the drug to be tested might be regarded as rather non-toxic too. Moreover, researchers would not create a version of an old drug that they would expect to be rather toxic but seek for such an improvement of an old one that yields, based on their prior knowledge, a non-toxic version. Therefore, again, there is a lack of convergence of higher pre-study probability of the hypothesis tested and the greater importance to avoid the error of falsely rejecting it. This pre-study belief that the physical hypothesis that the drug is non-toxic is rather true would remain valid even if one reformulated it as the alternative statistical hypothesis (and place the hypothesis of toxicity as the tested one). In such a case the more probable hypothesis would be the alternative, but then the importance of the IInd type of error would also change: it would become the error that is less important to avoid.

The negative influence of the above-discussed lack of convergence on $P_t$ can be numerically illustrated as follows. Assume, for example, that $P(H) = 0.6$, $\alpha = 0.05$ and $\beta = 0.1$, then $P_t$ equals 0.93. Now, should $\beta$ be recognized as more important to avoid at the cost of an increase of $\alpha$, for example, $\beta = 0.05$ and $\alpha = 0.1$, then $P_t$ would equal 0.72. The symmetric situation will take place when $P(H)$ is lower than one half and the error of the Ist type becomes more important to be avoided. The high discrepancy between $P(H)$ and the importance of error of falsely rejecting $H$ (e.g. low probability of the hypothesis and high importance of $\alpha$ compared to $\beta$) can make N-P unreliable epistemically. The example discussed in Section 5 suffices to show this: if $P(H) = 0.33$ and the more important error is $\alpha$, say, $\alpha = 0.05$ and $\beta = 0.94$,[24] then $P_t$ would equal 0.35 and in such a case, N-P's reliability will fail to meet the minimal level, hence the method would incline a researcher to a false assertion. The test would be at least minimally reliable only if $P(H)$ was sufficiently high under the assumption of the above-given error probabilities, or if the error of the IInd type was not that much less important compared to the error of the Ist type (under the assumption of the above-given hypothesis' probability).

## 8 Conclusions and Final Remarks

With the use of the concept of the predictive value, we checked whether N-P is at least minimally epistemically reliable for any possible and epistemically relevant situations of $\alpha$, $\beta$ and ratio $R$ ($P(H)$ alternatively), and we analyzed how pragmatic value-laden, asymmetric setting of error probabilities can influence this reliability and whether this influence can be controlled to have a non-negative impact. We found that N-P is incapable of warranting minimal epistemic reliability, except for the case of $R = 1$ ($P(H) = 0.5$). We also concluded that the character of influence of pragmatic value-laden asymmetry in magnitudes of error rates on the ER of N-P is accidental and impossible to be stipulated as principally positive.

The first finding of our study—certainty of at least minimal ER of N-P in every possible and epistemically relevant case only under the assumption of $R = 1$ is a potentially interesting philosophical problem for those who would expect a method of formulating scientific conclusions to lead more often to the truth than to false-hood by default. As we have shown if $R \neq 1$ warranty of at least minimal ER of N-P is impossible. Nevertheless, for users of N-P, this problem is not insurmountable because the method does not hinder the improvement of reliability by taking care of other aspects beyond the method itself. The level of $P_t$ depends especially on an experimental scheme (sample size) and the character of the studied phenomenon (its variance). A satisfactory reliability level is achievable by the use of the proper experimental scheme with a sufficient amount of data.

---

[24] Note that this represents the case of an underpowered, but, still, unbiased test.

The second finding of our study regarding the impact of pragmatic value-ladenness on N-P's ER has both potential philosophical and methodological consequences. It contributes to the topic of the role of values in science (see e.g. Elliott, Richards 2017). N-P captures the cognitive and non-cognitive factors that contribute to the research process, which has become the central topic of interest in post-Kuhnian theories of science. Within this approach, trends are developing that emphasize the fact that cognitive and social dynamics are inseparable elements of the cognitive act and the dynamics of scientific knowledge development (see Collins, Evans 2002; Kawalec 2020). An important example of such a direction is the so-called *Mode 2* science paradigm, which emphasizes social and economic factors in the formation of knowledge, in contrast with *Mode 1*, in which the formulation of scientific knowledge is motivated by the cognitive context alone (see e.g. Nowotny et al. 2001). N-P captures pragmatic factors as an element of the research process and by that, it is a classic example of how the general philosophical stance may be applied in precise methodological solutions.

We have shown how these pragmatic factors implemented in the form of the unequal setting of error probabilities may have a neutral, positive, or negative impact on the ER of N-P depending on the case of the physical hypothesis tested and the assumption about $R$. More importantly, we have shown that in the case of negative impact no methodological adjustment is available to neutralize it so in such cases the discussed pragmatic value-ladenness of N-P inevitably compromises the truth goal.

If epistemic and pragmatic aspects are assumed to be inseparable in the formation of knowledge and are implemented in statistical methodology in this type of tight relation, then a researcher has to challenge the fact that they inevitably are in mutual tension in some cases. Awareness of it is vital if a researcher wants to assess whether she wants to compromise or to give a higher rank to one of the two types of aspects and to communicate it to the society.

Even if a researcher intends to use N-P as a decision-theoretic tool and the general epistemic reliability indicator is not of primary concern, the outcomes of our study remain valid. In a decision-theoretic interpretation, they could be seen as an analysis of how often the method will, in the most general sense, lead to an erroneous decision, thus pragmatically not satisfactory in some way. Such an indication is certainly simplified, as it discerns no difference between the two types of pragmatic burden of an erroneous decision, but from what we have argued it follows that in N-P avoiding making extremely wrong decisions is in some cases inevitably conflicting with avoiding wrong decisions in general. This fact might, for example, serve as a justification of interpreting N-P as a type of minimax (minimization of worst possible loss) decision rule and by that to repel some doubts of whether N-P is a decision theory (see e.g. Szaniawski 1998) that rest on impossibility of calculating expectation of pragmatic loss.

**Funding Information** The author 1 gratefully acknowledges the support of the Polish National Science Center (Narodowe Centrum Nauki) under the grant no UMO-2015/17/N/HS1/02156. The author 2 gratefully acknowledges the support of the Minister of Science and Higher Education within the program under the name "Regional Initiative of Excellence" in 2019-2022, project number: 028/RID/2018/19, the amount of funding: 11 742 500 PLN.

## Declarations

**Conflict of interest** On behalf of all authors, the corresponding author states that there is no conflict of interest.